\documentclass[11pt,twoside,a4paper]{article}
\usepackage{amsfonts}
\usepackage{amssymb}
\usepackage{amsmath}

\begin{document}

\newcommand{\pa}{\partial}
\newcommand{\opa}{\overline\pa}
\newcommand{\ol}{\overline }

\numberwithin{equation}{section}

\newcommand\C{\mathbb{C}}  
\newcommand\R{\mathbb{R}}
\newcommand\Z{\mathbb{Z}}
\newcommand\N{\mathbb{N}}
\newcommand\PP{\mathbb{P}}

{\LARGE \centerline{Obstructions to finite dimensional cohomology}}
{\LARGE \centerline{of abstract Cauchy-Riemann complexes}}
\vspace{0.8cm}

\centerline{\textsc {Judith Brinkschulte\footnote{Universit\"at Leipzig, Mathematisches Institut, Augustusplatz 10, D-04109 Leipzig, Germany. 
E-mail: brinkschulte@math.uni-leipzig.de}}
 and C. Denson Hill
\footnote{Department of Mathematics, Stony Brook University, Stony Brook NY 11794, USA. E-mail: dhill@math.sunysb.edu\\
{\bf{Key words:}} cohomology of CR manifolds, non-solvability of tangential Cauchy-Riemann equations\\
{\bf{2010 Mathematics Subject Classification:}} 32V05, 32W10 }}

\vspace{0.5cm}

\begin{abstract} Let $M$ be a compact abstract $CR$ manifold of arbitrary $CR$ codimension. Under certain conditions on the Levi form we prove the infinite dimensionality of some global cohomology groups of $M$.
\end{abstract}

\vspace{0.5cm}

\section{Introduction}

Although there is a sizeable literature concerning various questions about $CR$ embeddable $CR$ manifolds, there appears to be very few results about the cohomology of abstract $CR$ manifolds. We consider a $\mathcal{C}^\infty$ smooth compact orientable abstract $CR$ manifold of type $(n,k)$.\\

Here an abstract $CR$ manifold of type $(n,k)$ is a triple $(M, HM, J)$, where $M$ is a smooth real manifold of dimension $2n+k$, $HM$ is a subbundle of rank $2n$ of the tangent bundle $TM$, and $J: HM \rightarrow HM$ is a smooth fiber preserving bundle isomorphism with $J^2= -\mathrm{Id}$. We also require that $J$ be formally integrable; i.e. that we have
$$\lbrack T^{0,1}M,T^{0,1}M\rbrack \subset T^{0,1}M$$
where 
$$ T^{0,1}M = \lbrace X+ iJX\mid X\in \Gamma(M,HM)\rbrace \subset \Gamma(M,\mathbb{C}TM),$$
with $\Gamma$ denoting smooth sections.

The $CR$ dimension of $M$ is $n\geq 1$ and the $CR$ codimension is $k\geq 1$.\\

We denote by $H^o M=\lbrace \xi\in T^\ast M\mid < X,\xi>=0, \forall X\in H_{\pi(\xi)}M\rbrace$ the {\it characteristic conormal bundle} of $M$. Here $\pi: T M \longrightarrow M$ is the natural projection. To each $\xi\in H^o_p M$, we associate the Levi form at $\xi:$
$$\mathcal{L}_p(\xi, X) = \xi(\lbrack J\tilde X, \tilde X\rbrack )= d\tilde\xi(X,JX) \ \mathrm{for} \ X\in H_p M$$
which is Hermitian for the complex structure of $H_p M$ defined by $J$. Here $\tilde \xi$ is a section of $H^o M$ extending $\xi$ and $\tilde X$ a section of $HM$ extending $X$. \\

We denote by $\opa_M$ the tangential Cauchy-Riemann operator on $M$. The associated cohomology groups of $\opa_M$ acting on smooth forms will be denoted by $H^{p,q}(M)$, $0\leq p\leq n+k,\ 0\leq q\leq n$. For more details on the $\opa_M$ complex, we refer the reader to \cite{HN1} or 
\cite{HN2}.\\

Our results are as follows.\\

\newtheorem{first}{Theorem}[section]
\begin{first}   \label{first}   \ \\
Let $M$ be a compact orientable abstract $CR$ manifold of type $(n,k)$. Assume that there exists a point $p_0\in M$ and a characteristic conormal direction $\xi\in H^o_{p_0}M$ such that the Levi form $\mathcal{L}_{p_0}(\xi,\cdot)$ has $q$ negative and $n-q$ positive eigenvalues. Then for $0\leq p\leq n+k$, the following holds: Either $H^{p,q}(M)$ or $H^{p,q+1}(M)$ is infinite dimensional, and either $H^{p,n-q}(M)$ or $H^{p,n-q+1}(M)$ is infinite dimensional.
\end{first}

The theorem is proved in section 2. Although here we are proving the infinite dimensionality
of certain (global) cohomology groups of $M$, our argument follows the pattern of M. Nacinovich \cite{N}, where the emphasis was on demonstrating the absence of the (local) Poincar\'e lemma.\\

The following two theorems are consequences of Theorem \ref{first}. The simple arguments are given at the end of section 2.

\newtheorem{second}[first]{Theorem}
\begin{second}   \label{second}   \ \\
Let $M$ be a compact orientable abstract $CR$ manifold of type $(n,1)$. Assume that at each point $x\in M$, there exists a characteristic conormal direction $\xi\in H_x^o(M)$ such that $\mathcal{L}_x(\xi,\cdot)$ has $q$ negative and $n-q$ positive eigenvalues. If moreover $2q\not= n-1$, then $H^{p,q}(M)$ is infinite dimensional; and if $2q\not= n+1$, then $H^{p,n-q}(M)$ is infinite dimensional, $0\leq p\leq n+1$.
\end{second}

\newtheorem{third}[first]{Theorem}
\begin{third}   \label{third}   \ \\
Let $M$ be a compact orientable abstract $CR$ manifold of type $(n,k)$ with $n$ even. For $q=\frac{n}{2}$ we assume that  at each point $x\in M$ and every characteristic conormal direction $\xi\in H_x^o(M)\setminus\lbrace 0\rbrace$ the Levi form $\mathcal{L}_x(\xi,\cdot)$ has $q$ negative and $q$ positive eigenvalues. Then $H^{p,q}(M)$ is infinite dimensional, $0\leq p\leq n+k$.
\end{third}

In the case where $M$ is $CR$ embedded in some ambient complex manifold, related local and global results have been discussed in \cite{AH1}, \cite{AH2}, \cite{AFN} and \cite{HN2}.

\section{Proofs of the theorems}

Our proof of Theorem \ref{first} relies on a well known construction for $CR$ embedded $CR$ manifolds at a point where there exists a characteristic conormal direction such that the associated Levi form has exactly $q$ negative and $n-q$ positive eigenvalues. For the reader's convenience, we now sketch this construction in the case of a hypersurface in $\mathbb{C}^{n+1}$. \\

So let $S\ni 0$ be a piece of a smooth real hypersurface in $\C^{n+1}$ such that $\mathcal{L}_0(\xi,\cdot)$ has $q$ negative and $n-q$ positive eigenvalues for some characteristic conormal direction $\xi$. Then we can choose a local real defining function $\rho$ of $S$ of the form
$$\rho(z)= \mathrm{Im}(z_{n+1}) - h(z) \quad \mathrm{with}\ h(z) = O(\vert z\vert^2).$$
Here $O(\vert z\vert^\ell)$ denotes a term vanishing to order $\ell$ at the point $z=0$.
Moreover, after a holomorphic change of coordinates, we may assume
$$h(z) = \sum_{\alpha=1}^q \vert z_\alpha\vert^2 - \sum_{\alpha =q+1}^n\vert z_\alpha\vert^2 + O(\vert z\vert^3) \quad \mathrm{at}\ 0.$$
Set 
$$\phi(z)= -i\mathrm{Re}(z_{n+1}) + h(z) - 2 \sum_{\alpha=1}^q \vert z_\alpha\vert^2 -  (\mathrm{Re}(z_{n+1}) +ih(z))^2.$$

Then
$$\mathrm{Re}\phi(z) \leq -\frac{1}{2} (\sum_{\alpha=1}^n\vert z_\alpha\vert^2 +(\mathrm{Re}(z_{n+1}))^2)\quad \mathrm{near}\ 0.$$
For $\lambda > 0$ we then define the following "peak forms"
$$f_\lambda = e^{\lambda\phi} dz_1\wedge\ldots\wedge dz_p\wedge d\ol z_1\wedge\ldots d\ol z_q,$$
which defines a smooth $(p,q)$-form on $S$ satisfying $\opa_S f_\lambda =0$ (note that $\mathrm{
Re}(z_{n+1}) +ih$ is the restriction to $S$ of the holomorphic function $z_{n+1}$).\\

Similarly we set
$$\psi (z) = i\mathrm{Re}(z_{n+1})-h(z) -2\sum_{\alpha=q+1}^n \vert z_\alpha\vert^2 -(\mathrm{Re}(z_{n+1}) +ih(z))^2.$$
Then we also have
$$\mathrm{Re}\psi(z) \leq -\frac{1}{2}(\sum_{\alpha=1}^n\vert z_\alpha\vert^2+ (\mathrm{Re}(z_{n+1}))^2)\quad \mathrm{near}\ 0,$$
and we define the following "peak forms" of degree $(n+1-p,n-q)$ on $S$:
$$g_\lambda = e^{\lambda\psi} dz_{p+1}\wedge\ldots \wedge dz_{n+1} \wedge d\ol z_{q+1}\wedge\ldots\wedge d\ol z_n.$$
Again we have $\opa_S g_\lambda =0.$\\

In the proof of the nonvalidity of the Poincar\'e lemma for the $\opa_S$-operator, the forms $f_\lambda$ and $g_\lambda$ play an essential role, because their properties contradict the existence of certain a priori estimates. Also our proof of Theorem \ref{first} is based on the existence of forms with the analogous properties up to some terms vanishing to infinite order at the point under consideration.\\

For more details on the construction of the corresponding functions and forms in the higher codimensional situation, we refer the reader to the paper $\cite[\mathrm{pp.} 388 ff.]{AFN}$. \\

{\it Proof of Theorem \ref{first}.}  

Let us first  consider the case $q> 0$, and  assume by contradiction that both $H^{p,q}(M)$ and $H^{p,q+1}(M)$ are finite dimensional. In order to make the proof as clear as possible, we first assume that $k=1$ ($CR$ manifold of hypersurface type).\\

By the assumption that $H^{p,q}(M)$ is finite dimensional we get that
$$\opa_M : \mathcal{C}^\infty_{p,q-1}(M) \longrightarrow \mathcal{C}^\infty_{p,q}(M) $$
has closed range. Then Banach's open mapping theorem implies that there exist a constant $C_1> 0$ and an integer $m_1 >0$ such that for all $f\in\opa_M\mathcal{C}^\infty_{p,q-1}(M) $ there exists $u\in \mathcal{C}^\infty_{p,q-1}(M) $ satisfying $\opa_M u=f$ and 
\begin{equation}  \label{apriori1}
\Vert u\Vert_0 \leq C_1 \Vert f\Vert_{m_1}.
\end{equation}
Here $\Vert \ \Vert_m$ denotes the usual $\mathcal{C}^m$-norm on $\mathcal{C}^\infty_{\cdot,\cdot}(M)$.\\

Reasoning as before, the assumption that $H^{p,q+1}(M)$ is finite dimensional implies that there exist a constant $C_2 > 0$ and and integer $m_2 > 0$ such that for all $g\in\opa_M\mathcal{C}^\infty_{p,q}(M) $ there exists $h\in \mathcal{C}^\infty_{p,q}(M) $ satisfying $\opa_M h=g$ and 
\begin{equation}  \label{apriori2}
\Vert h\Vert_{m_1} \leq C_2 \Vert g\Vert_{m_2}.
\end{equation}

Using Stokes' formula, we have for every $f=\opa_M u\in \opa_M\mathcal{C}^\infty_{p,q-1}(M)$ and every $g\in \mathcal{C}^\infty_{n+1-p,n-q}(M)$:
$$\int_M f\wedge g = \int_M \opa_M u\wedge g = (-1)^{p+q}\int_M u\wedge\opa_M g.$$

Hence (\ref{apriori1}) implies 
\begin{equation}  \label{apriori3}
\vert\int_M f\wedge g\vert \lesssim C_1 \Vert f\Vert_{m_1} \cdot \Vert \opa_M g\Vert_0
\end{equation}
for every $f\in  \opa_M\mathcal{C}^\infty_{p,q-1}(M)$ and every $g\in \mathcal{C}^\infty_{n+1-p,n-q}(M)$. Here $a\lesssim b$ means that there exists a constant $C > 0$ such that $a \leq C\cdot b$. \\

Now let $l:= \dim H^{p,q}(M) < + \infty$, and let $\Omega$ be an open neighborhood of $p_0\in M$ such that for every point $x\in \Omega$, there exists a characteristic conormal direction $\xi_x$ such that $\mathcal{L}_x(\xi_x,\cdot)$ has $q$ negative and $n-q$ positive eigenvalues.\\

We choose $l$ different points $p_1,\ldots, p_l$ inside $\Omega$, all different from $p_0$. Moreover, we choose cut-off functions $\chi_j,\ j=0,\ldots, l$, with $\chi_j\equiv 1$ near $p_j$, such that the $\chi_j$'s have disjoint supports. For each $0\leq j\leq l$, we then make the following construction:\\

Choose local coordinates $z_1=x_1+ix_{n+1},\ldots,z_n= x_n+i x_{2n}, x_{2n+1}$ for $M$ so that $p_j$ becomes the origin. By the formal Cauchy-Kowalewski procedure, we can find smooth complex valued functions $\varphi = (\varphi_1,\ldots, \varphi_{n+1})$ in an open neighborhood $U$ of $0$ with each $\varphi_i(0)=0,\ d\varphi_1\wedge\ldots\wedge\varphi_{n+1}\not= 0$ in $U$, and such that $\opa_M \varphi_i$ vanishes to infinite order at $0$. Then $\varphi: U\longrightarrow \mathbb{C}^{n+1}$ gives a smooth local embedding $\tilde M = \varphi(U)$ of $M$ into $\mathbb{C}^{n+1}$. On $\tilde M$ there is the $CR$ structure induced from $\mathbb{C}^{n+1}$; it agrees to infinite order at $0$ with the original $CR$ structure on $M$. In particular $\tilde M$ is a smooth real hypersurface in $\mathbb{C}^{n+1}$ which is strictly $q$-convex and strictly $(n-q)$-concave with respect to the induced $CR$ structure. As explained in the paragraphs preceeding this proof  this means that after possibly shrinking $U$, there are smooth complex valued functions $\phi_j$ and $\psi_j$ on $U$ with $\opa_M\phi_j$ and $\opa_M\psi_j$ vanishing to infinite order at $0$ satisfying
\begin{equation}  \label{phi}
\mathrm{Re}\phi_j \leq -\frac{1}{2}\vert x\vert^2 \quad\mathrm{on} \ U,
\end{equation}

\begin{equation}  \label{psi}
\mathrm{Re}\psi_j \leq -\frac{1}{2}\vert x\vert^2 \quad\mathrm{on} \ U
\end{equation}
and
\begin{equation}  \label{phi+psi}
\phi_j +\psi_j = -2\vert x\vert^2 + O(\vert x\vert^3)
\end{equation}
(one constructs the corresponding functions $\phi$ and $\psi$ on $\tilde M$ and considers the pull-back under $\varphi$.)\\

Moreover, $T^\ast M$ is spanned by forms
$$\omega_1 = dz_1 + O(\vert x\vert^\infty),\ldots,\omega_n = dz_n + O(\vert x\vert^\infty),$$
$$\ol\omega_1= d\ol z_1 + O(\vert x\vert^\infty),\ldots,\ol\omega_n = d\ol z_n = O(\vert x\vert^\infty),$$
$$ \theta = dx_{2n+1} + O(\vert x\vert^\infty)$$
 which are $d$-closed to infinite order at $0$ (here, of course, $\Lambda T^{0,1}M$ is spanned by $\ol\omega_1,\ldots,\ol\omega_n$). Following again $\cite{AFN}$ or $\cite{HN2}$, by the geometric condition on the Levi-form of $M$ we may also assume that $\opa_M\phi_j \wedge\ol\omega_1\wedge\ldots\wedge\ol\omega_q$ and $\opa_M\psi_j\wedge\ol\omega_{q+1}\wedge\ldots\wedge\ol\omega_n$ vanish to infinite order at $0$.\\

For each real $\lambda > 0$ we now define 
$$ f_j^\lambda = \chi_j e^{\lambda\phi_j} \omega_1\wedge\ldots\wedge\omega_p\wedge\ol\omega_1\ldots\wedge\ol\omega_q.$$
This is a smooth $(p,q)$-form on $M$. Moreover the properties of $\phi_j$ imply
that $\opa_M(f_j^\lambda)$ is rapidly decreasing with respect to $\lambda$ in the topology of $\mathcal{C}^\infty(M)$ as $\lambda$ tends to infinity. Indeed, by (\ref{phi}) the function $\exp(\lambda\phi_j)$ , and any derivative of it with respect to $x$, is rapidly decreasing as $\lambda\rightarrow +\infty$, while all other terms, and their derivatives with respect to $x$, have only polynomial growth in $\lambda$.\\

We also set
$$ g_j^\lambda = \chi_j e^{\lambda\psi_j} \omega_{p+1} \wedge\ldots\wedge\omega_n\wedge\theta\wedge\ol\omega_{q+1}\wedge\ldots\wedge\ol\omega_n.$$
Then, arguing as before, $\opa_M(g_j^\lambda)$ is rapidly decreasing with respect to $\lambda$ in the topology of $\mathcal{C}^\infty(M)$ as $\lambda$ tends to infinity.\\

Next, we solve $\opa_M u_j^\lambda = \opa_M f_j^\lambda$ with an estimate
\begin{equation}  \label{apriori4}
\Vert u_j^\lambda\Vert_{m_1} \leq C_2 \Vert \opa_M f_j^\lambda\Vert_{m_2} ,
\end{equation}
using (\ref{apriori2}). Hence 
$\Vert u_j^\lambda\Vert_{m_1}$ is rapidly decreasing with respect to $\lambda$.
 Defining $\tilde f_j^\lambda = f_j^\lambda - u_j^\lambda$, we obtain a smooth, $\opa_M$-closed $(p,q)$-form on $M$. \\
 
 Since $\mathrm{dim}H^{p,q}(M)=l$, there exist constants $c_0^\lambda,\ldots, c_l^\lambda$, not all equal to zero, such that
 $$c_0^\lambda\tilde f_0^\lambda + \ldots + c_l^\lambda \tilde f_l^\lambda \in\mathrm{Im}\opa_M.$$
 To get a contradiction, we are going to use the estimate (\ref{apriori3}) with $f=\sum_{j=0}^l c_j^\lambda \tilde f_j^\lambda$ and $g= \sum_{j=0}^l \ol c_j^\lambda g_j^\lambda$. We have
 
 \begin{eqnarray}
 \int_M f\wedge g \qquad\qquad & = &  \label{fwedgeg} \\
 \int_M (\sum_{j=0}^l c_j^\lambda \tilde f_j^\lambda)\wedge (\sum_{j=0}^l \ol c_j^\lambda g_j^\lambda) & = & \int_M (\sum_{j=0}^l c_j^\lambda ( f_j^\lambda-u_j^\lambda))\wedge (\sum_{j=0}^l \ol c_j^\lambda g_j^\lambda)\nonumber\\
& = & \sum_{j=0}^l\vert c_j^\lambda\vert^2 \int_M f_j^\lambda\wedge g_j^\lambda - \int_M \sum_{i,j=0}^l c_i^\lambda \ol c_j^\lambda u_i^\lambda\wedge g_j^\lambda\nonumber
\end{eqnarray}
Note that for the third equality, we have used that the $\chi_j$'s have disjoint supports.\\

We are now going to estimate the term on the right of (\ref{fwedgeg}). We have
\begin{eqnarray*}
\int_M f_j^\lambda\wedge g_j^\lambda & = & \int_M \chi_j^2 e^{\lambda(\phi_j + \psi_j)} \omega_1\wedge\ldots\wedge\omega_n\wedge\theta\wedge\ol\omega_1\wedge\ldots\wedge\ol\omega_n \\
& = & \int_M \lbrace \chi_j^2 e^{\lambda(-2\vert x\vert^2 + O(\vert x\vert^3)}  + O(\vert x\vert) \rbrace dz_1\wedge\ldots\wedge dz_n\wedge d\ol z_1\wedge\ldots\wedge d\ol z_n\wedge dx_{2n+1}.
\end{eqnarray*}
Making the change of variables $y=\sqrt\lambda x$, and afterwards changing the name of $y$ back to $x$,  we get
\begin{equation} \nonumber
\int_M f_j^\lambda\wedge g_j^\lambda = \lambda^{-n-\frac{1}{2}} \lbrace \int_M \chi_j^2(\frac{x}{\sqrt\lambda}) e^{-2\vert x\vert^2 + O(\lambda^{-\frac{1}{2}})}dz_1\wedge\ldots\wedge dz_n\wedge d\ol z_1\wedge\ldots\wedge d\ol z_n\wedge dx_{2n+1} + O(\lambda^{-\frac{1}{2}}) \rbrace
.
\end{equation}

Therefore we obtain
\begin{equation} \label{firstintegral}
\vert\int_M f_j^\lambda\wedge g_j^\lambda\vert \geq c \lambda^{-n-\frac{1}{2}}
\end{equation}
for some constant $c > 0$.\\

Also we can use (\ref{apriori4}) to get
\begin{eqnarray*}
\vert \int_M \sum_{i,j=0}^l c_i^\lambda \ol c_j^\lambda u_i^\lambda\wedge g_j^\lambda \vert & \lesssim & 
\sum_{j=0}^l\vert c_j^\lambda\vert^2 \sup_{i,j} (\Vert u_i^\lambda\Vert_0\cdot\Vert g_j^\lambda\Vert_0)\\
& \lesssim & \sum_{j=0}^l\vert c_j^\lambda\vert^2 \sup_{i,j} (\Vert \opa_M f_i^\lambda\Vert_{m_2}\cdot\Vert g_j^\lambda\Vert_0).
\end{eqnarray*}

Now $\Vert\opa_M f_i^\lambda\Vert_{m_2}$ is rapidly decreasing with respect to $\lambda$, whereas $\Vert g_j^\lambda\Vert_0$ is of polynomial growth with respect to $\lambda$, hence we get 
\begin{equation}  \nonumber
\vert \int_M \sum_{i,j=0}^l c_i^\lambda \ol c_j^\lambda u_i^\lambda\wedge g_j^\lambda \vert \leq 
\sum_{j=0}^l\vert c_j^\lambda\vert^2 \lambda^{-n-1}
\end{equation}
for sufficiently large $\lambda$. Combining this with (\ref{firstintegral}), we get
\begin{equation} \label{integral}
\vert\int_M f\wedge g\vert \geq \frac{c}{2}\sum_{j=0}^l\vert c_j^\lambda\vert^2  \lambda^{-n-\frac{1}{2}}
\end{equation}  
for sufficiently large $\lambda$.

On the other hand,
using (\ref{apriori3}), we can estimate $\int_M f\wedge g$ as follows:

$$\vert \int_M f\wedge g\vert  \leq  C_1 \Vert f\Vert_{m_1} \cdot \Vert \opa_M g\Vert_0$$
$$ \lesssim  \sum_{j=0}^l \vert c_j^\lambda\vert^2 \sup_{i,j} (\Vert\tilde f_j^\lambda\Vert_{m_1} \cdot \Vert \opa_M g_j^\lambda\Vert_0 ) $$
$$ \lesssim  \sum_{j=0}^l \vert c_j^\lambda\vert^2 \sup_{i,j} (\Vert  f_j^\lambda\Vert_{m_2+1} \cdot \Vert \opa_M g_j^\lambda\Vert_0).$$

Since $\Vert f_j^\lambda\Vert_{m_2+1}$ is of polynomial growth in $\lambda$ whereas $\Vert \opa_M g_j^\lambda\Vert_0$ is rapidly decreasing with respect to $\lambda$, we get that
$$\vert \int_M f\wedge g\vert  \lesssim \sum_{j=0}^l \vert c_j^\lambda\vert^2 \lambda^{-n-1}.$$
This contradicts (\ref{integral}) and therefore proves that either $H^{p,q}(M)$ or $H^{p,q+1}(M)$ has to be infinite dimensional. \\

Now, replacing $\xi$ by $-\xi$,  it also follows that either $H^{p,n-q}(M)$ or $H^{p,n-q+1}(M)$ is infinite dimensional.\\

For $q=0$, the statement is essentially Boutet de Monvel's result \cite{BdM}: In this case, $M$ is strictly pseudoconvex at $p_0$. If $H^{p,1}(M)$ was finite dimensional, then in particular the range of $\opa_M$ was closed in $\mathcal{C}^\infty_{p,1}(M)$. But then one can construct infnitely many linearly independent $CR$ functions on $M$ as in \cite{BdM}.\\

Also, the Levi-form $\mathcal{L}_{p_o}(-\xi,\cdot)$ has $n> 0$ negative and $0$ positive eigenvalues. By what already proved, we therefore know that $H^{p,n}(M)$ is infinite dimensional (note that $H^{p,n+1}(M)$ is always zero).\\

If $k > 1$,  the proof is essentially as before, with $\mathbb{C}^{n+1}$ replaced by $\mathbb{C}^{ n+k}$. The crucial point is to observe that the approximate $CR$ embedding $\tilde M$ in $\mathbb{C}^{n+k}$, which now has real codimension $k$, is contained in a hypersurface which is strictly $q$-convex and strictly $(n-q)$-concave. 
$\xi$ then corresponds to + or -- the conormal to the hypersurface at $p_0$.
As before, this gives us the existence of smooth functions $\phi_j$ and $\psi_j$ with the same properties that were essential in the proof for $k=1$.
\hfill$\square$\\

{\it Proof of Theorem \ref{second}.} 

The theorem follows immediately from Theorem \ref{first}. Indeed,  it suffices to note that the assumptions on $M$ imply that the classical conditions $Y(q+1)$ and $Y(n-q+1)$ are satisfied. Hence the $\opa_M$-complex is $\frac{1}{2}$-subelliptic in degree $(p,q+1)$ and $(p,n-q+1)$ (see $\cite[\mathrm{Theorem}\ 5.4.9]{FK}$), hence $H^{p,q+1}(M)$ and $H^{p,n-q+1}(M)$ are finite dimensional. \hfill$\square$\\

{\it Proof of Theorem \ref{third}.} 
The assumptions on $M$ imply that $M$ is $q$-pseudoconcave (see \cite{HN1} for the definition), hence the $\opa_M$-complex is $\epsilon$-subelliptic in degree $(p,q+1)$  for some $\epsilon > 0$ (see \cite{HN1} for the proof), hence $H^{p,q+1}(M)$ is finite dimensional. Again the statement then follows from Theorem \ref{first}. \hfill$\square$\\

\section{Corollaries and remarks}

We would like to emphasize that Theorems \ref{first}, \ref{second}, \ref{third} are valid for an {\it abstract} $CR$ manifold $M$, which might possibly be not even locally $CR$ embeddable at any point. However, it is of some interest to consider the situation where $M$ is globally $CR$ embedded as a generic $CR$ submanifold of some complex manifold $X$, and ask what these theorems imply about the pair $(X,M)$. Then the complex dimension of $X$ is $n+k$, and we have the usual Dolbeault cohomology groups $H^{p,q}(X)$, as well as the Dolbeault-like cohomology groups $H^{p,q}(X,\mathcal{I})$. The latter consists of smooth $\opa$-closed forms on $X$ modulo smooth $\opa$-exact forms on $X$, in which all forms are required to have zero Cauchy data along the submanifold $M$. (Think of the real codimension $k$ of $M$ in $X$ as corresponding to $k$ "time variables".) More precisely, we consider the sheaf $\mathcal{I}_M$ of germs of $\mathcal{C}^\infty$ functions on $X$ which vanish on $M$. Then we denote by $\mathcal{I}$ the sheaf of $\mathcal{C}^\infty_{\cdot,\cdot}(X)$-modules which is locally generated by $\mathcal{I}_M$ and $\opa\mathcal{I}_M$.\\

The interpretation of $H^{p,q+1}(X,\mathcal{I})$ is that it is the obstruction to the solvability of the general inhomogeneous Cauchy problem
\begin{equation}  \label{cauchy}
\left\lbrace \begin{array}{c}
\opa u = f\ \mathrm{on}\ X\\
u=u_0\ \mathrm{on}\ M. 
\end{array}  \right.
\end{equation}
Here $f$ is a given smooth $\opa$-closed $(p,q+1)$-form on $X$, $u_0$ is a given smooth $\opa_M$-closed tangential $(p,q)$-form on $M$, and it is assumed that the data $\lbrace f, u_0\rbrace$ are compatible (see \cite[p. 350--351]{AH1}). The desired solution $u$ to the problem (\ref{cauchy}) should be a smooth $(p,q)$-form on $X$. Then the solvability of (\ref{cauchy}) for all compatible data is equivalent to the vanishing of $H^{p,q+1}(X,\mathcal{I})$. If, for example, $H^{p,q+1}(X,\mathcal{I})$ is infinite dimensional, it means that there is an infinite dimensional set of equivalent classes of data for which (\ref{cauchy}) has no solution (see \cite{AH1}). 
 From Theorem \ref{first} and standard exact sequences, such as the Mayer-Vietoris sequence, we obtain the following Corollary:

\newtheorem{cor}{Corollary}[section]
\begin{cor}   \label{cor}   \ \\
With $M$ as in Theorem \ref{first}, assume $X$ is either compact or Stein. Then
\indent (a) either $H^{p,q+1}(X,\mathcal{I})$ or $H^{p,q+2}(X,\mathcal{I})$ is infinite dimensional.\\
And\\
\indent (b) either $H^{p,n-q+1}(X,\mathcal{I})$ or $H^{p,n-q+2}(X,\mathcal{I})$ is infinite dimensional.
\end{cor}

{\it Proof.} If $X$ is compact, then we have $\mathrm{dim} H^{r,s}(X)< + \infty$ for $0\leq r,s\leq n+k$, whereas $H^{r,s}(X)=0$ for $0\leq r\leq n+k, 1\leq s\leq n+k$ if $X$ is Stein.
Therefore we may use the following long exact sequence

$$\ldots \rightarrow H^{r,s}(X) \rightarrow H^{r,s}(M)\rightarrow H^{r,s+1}(X,\mathcal{I})                           \rightarrow H^{r,s+1}(X)\rightarrow\ldots $$

and Theorem \ref{first} to conclude. \hfill$\square$\\

 In the special case $k=1$, we may assume $M$ divides $X$ into complex manifolds-with-boundary, which we call $X^-$ and $X^+$ (the common boundary is, of course, $M$). Then, roughly speaking, the cohomology group $H^{p,j+1}(X,\mathcal{I})$ is isomorphic to the direct sum of $H^{p,j}(X^-)$ and $H^{p,j}(X^+)$ modulo some global Dolbeault cohomology groups of $X$. Therefore the hypersurface case of Corollary \ref{cor} then breaks down into:\\
 
 \noindent (a) at least one of $H^{p,q}(X^+)$, $H^{p,q}(X^-)$, $H^{p,q+1}(X^+)$, $H^{p,q+1}(X^-)$ is infinite dimensional.\\
 And\\
 (b) at least one of $H^{p,n-q}(X^+)$, $H^{p,n-q}(X^-)$, $H^{p,n-q+1}(X^+)$, $H^{p,n-q+1}(X^-)$ is infinite dimensional.\\

 We should emphasize here that the above corollary requires  a hypothesis on $M$  at only {\it one} single point $p_0$ on $M$ and  in only {\it one} single characteristic conormal direction $\xi$. This has the following consequence: Suppose $M$ is generically $CR$ embedded in $X$, as above, with $X$ either compact or Stein, but that initially no other hypotheses are made about $M$. Then the situation is initially whatever it is. But if now we make arbitrarily small smooth modifications of $M$ at a few points, we can produce a modified $CR$ manifold $\tilde M$, such that for the new pair $(X,\tilde M)$, there is a plethora of infinite dimensional cohomology groups $H^{p,\ast}(X,\tilde{\mathcal{I}})$.\\
 
 Note that in Theorems \ref{second} and \ref{third} the situation is quite different. In those theorems some hypothesis is needed at each point of $M$, which puts a much greater constraint of the "shape" of $M$, and we are then in a territory that is less novel and has been much more discussed in the literature.\\

Indeed,  with $M$ as in Theorem \ref{second}, and $X$ compact, we know from \cite[p. 805]{AH2} that it is $H^{p,q}(X^-)$ and $H^{p,n-q}(X^+)$ that are infinite dimensional.
 Also $H^{p,j}(X^-)$ is finite for $j\not= q$ and $H^{p,j}(X^+)$ is finite for $j\not= n-q$. So in this context, given the finiteness theorems proved in \cite{AH2}, what Theorem \ref{second} provides us in most cases is just a new proof of the infinite dimensionality of $H^{p,q}(X^-)$ and $H^{p,n-q}(X^+)$. We should also recall from \cite{AH2} that when $2q\not= n$, we therefore have a good one sided global Cauchy problem in degree $q$ from the $X^-$ side, and another one in degree $n-q$ from the $X^+$ side. Both these Cauchy problems are almost always solvable. If $2q= n$, then we have an almost always solvable Riemann-Hilbert problem in degree $q=n-q$.\\

Now with $M$ as in Theorem \ref{third} and $X$  compact, $M$ is a maximally pseudoconcave generic $CR$ submanifold of $X$, of codimension $k$. Theorem \ref{third} gives us a new proof (in the maximally pseudoconcave case) of the infinite dimensionality of $H^{p,q}(M)$, which is related to Theorem 4.2 in \cite{HN1}. In that situation the global solvability of the inhomogeneous Cauchy problem (\ref{cauchy}) is obstructed by the infinite dimensional $H^{p,q+1}(X,\mathcal{I})$.\\

\section{Examples} 

Standard examples of compact hypersurfaces satisfying the assumptions of Theorem \ref{second} are the real projective hypersurfaces
$$M = \lbrace (z_0 : z_1 : \ldots : z_{n+1})\in\C\PP^{n+1} \mid \mathrm{Im}(z_0 z_{n+1}) = \vert z_1\vert^1 + \ldots + \vert z_q\vert^2 - \vert z_{q+1}\vert^2 - \ldots - \vert z_n\vert^2 \rbrace.$$
Various other examples of $CR$ manifolds satisfying the assumptions of Theorem \ref{second} or Theorem \ref{third} have been constructed by C. Medori and M. Nacinovich: They continued the investigations of Tanaka and developed the algebraic theory of Levi-Tanaka algebras in order to construct homogeneous $CR$ manifolds of arbitrary codimension $k$. In \cite[$\mathrm{Theorem}\ 4.5$]{MN1} they showed that if the Levi-Tanaka algebra $\mathfrak{g}$ is semisimple, then the associated homogeneous $CR$ manifold $M_{\mathfrak{g}}$ is compact.
Moreover, in \cite{MN2} it was proved that the Levi-form of $M$ is nondegenerate if and only if the Levi-Tanaka algebra is finite dimensional. A complete classification of semisimple Levi-Tanaka algebras was also given in \cite{MN2}. 
So in those examples, we get from Theorems \ref{second} and \ref{third} that for fixed $p$, one or two cohomology groups are infinite dimensional, while all others are finite dimensional.\\

Here is a method how to construct compact $CR$ manifolds with at least approximately half of its cohomology groups being infinite dimensional:\\

We start with a compact $CR$ manifold $M$ of arbitrary type $(n,k)$, which we assume to be generically $CR$  embedded into some complex manifold $X$. Now we cut out a small piece of $M$ near a given point and glue in a small modification, arranging that $q=0$ at one point, that $q=1$ at another point, ..., and that $q=n$ at still another point (all happening locally in $\C^{n+k}$). We denote the modified $CR$ manifold by $\tilde M$. Then, using Theorem \ref{first}, we obtain that either $H^{p,j}(\tilde M)$ or $H^{p,j+1}(\tilde M)$ must be infinite dimensional for all $j=0,\ldots ,n$.\\

It is, however, far more difficult to produce examples of {\it abstract} $CR$ structures having a characteristic conormal direction whose associated Levi-form is nondegenerate. Examples exist, but they are few. Here we would like to mention the following example from \cite[Theorem \ 6.16]{HN3}:

Let
$Q\subset\C\PP^{n+1},\ n\geq 1,$ be the hyperquadric defined by
 $$Q = \lbrace z_0\ol z_0 + z_1\ol z_1 = z_2\ol z_2 + \ldots + z_{n+1}\ol z_{n+1}\rbrace.$$
Then one can find a new $CR$ structure on $Q$, which is not locally $CR$-embeddable at all points of the divisor $D=\lbrace z_0 =0\rbrace$. We denote $Q$ with this new $CR$ strucure by $\tilde Q$. The $CR$ structure on $\tilde Q$ is such that at each point $x\in M$, there exists a characteristic  conormal direction such that $\mathcal{L}_x(\xi,\cdot)$ has $1$ negative and $n-1$ positive eigenvalues, i.e. $\tilde Q$ satisfies the assumptions of Theorem \ref{second} with $q=1$. In this situation, Theorem \ref{second} yields that if $n\not= 3$, then $H^{p,1}(\tilde Q)$ is infinite dimensional, and if $n\not= 1$, then $H^{p,n-1}(\tilde Q)$ is infinite dimensional, $0\leq p\leq n+1$. This is a new result.\\

\end{document}